\documentclass[11pt,a4paper]{article}

\usepackage[utf8]{inputenc}
\usepackage[T1]{fontenc}

\usepackage[margin=1in,includefoot]{geometry}

\usepackage{amsmath,amssymb,amsthm}
\newtheorem{theorem}{Theorem}[section]

\usepackage{graphicx}
\usepackage{subfig}

\usepackage{enumitem}

\usepackage[colorlinks=true,
            linkcolor=blue,
            urlcolor=blue,
            citecolor=blue]{hyperref}

\usepackage{authblk}

\newcommand{\E}{\mathbb{E}}
\newcommand{\Prob}{\mathbb{P}}
\newcommand{\ov}[1]{\overline{#1}}

\theoremstyle{definition}
\newtheorem{defn}[theorem]{Definition}
\newtheorem{ex}[theorem]{Example}
\newtheorem{thm}[theorem]{Theorem}
\newtheorem{cor}[theorem]{Corollary}

\theoremstyle{remark}
\newtheorem{rem}[theorem]{Remark}

\newcommand{\keywords}[1]{\par\noindent\textbf{Keywords:} #1}

\title{Any random variable with right-unbounded \\
 distributional support is the minimum of \\
 independent and very heavy-tailed random \\
 variables} 

\author[1,2]{Sergey Foss\thanks{sergueiorfoss25@gmail.com}}
\author[2]{Anton Tarasenko\thanks{tarasenko@math.nsc.ru}}
\author[3,4]{Georgiy Krivtsov\thanks{gkrivtsov@nes.ru}}

\affil[1]{Heriot-Watt University}
\affil[2]{Sobolev Institute of Mathematics}
\affil[3]{Novosibirsk State University}
\affil[4]{New Economic School}

\date{}

\begin{document}
\maketitle

\begin{abstract}
A random variable $\xi$ has a {\it light-tailed} distribution (for short: is light-tailed) if it possesses a finite exponential moment, $\E \exp (\lambda \xi) <\infty$ for some $\lambda >0$, and has a {\it heavy-tailed} distribution (is heavy-tailed) if 
$\E \exp (\lambda\xi) = \infty$, for all $\lambda>0$.
In (Leipus et al., AIMS Mathematics, 2023), the authors presented a particular example of a light-tailed random variable that is the minimum of two independent heavy-tailed random variables.
We will show that this phenomenon is universal: {\it any} light-tailed  random variable with right-unbounded support may be represented as the minimum of two independent heavy-tailed random variables. Moreover, a more general fact holds: these two independent random variables may have as heavy-tailed distributions as one wishes. 
Further, we will extend the latter result onto the minimum of any finite number of independent random variables.
We will also comment on possible generalizations of our result to the case of dependent random variables.
\end{abstract}


\section{Introduction and results} 


Questions about the closure of certain distribution classes under various operations (summation, taking maximum or minimum, etc.) arise naturally in the tail asymptotics analysis and in a variety of applications, see, e.g. \cite{EKM}---\cite{XFW} and the references therein.

We mostly restrict our attention to independent random variables, and 
there are several known simple closure properties. 
For example, the sum of two variables is heavy-tailed if at least one of the summands is heavy-tailed, and light-tailed otherwise; the sum is long-tailed if at least one of the  summands is long-tailed; the minimum and the maximum of light-tailed random variables are light-tailed; the minimum and the maximum of long-tailed random variables are long-tailed; and the maximum of heavy-tailed random variables is heavy-tailed. 

However, it was shown by example in \cite{LSK1} that the minimum of two independent heavy-tailed random variables may be light-tailed. 
The example, however, is somewhat involved and does not fully elucidate
the mechanism behind this phenomenon.
This was our motivation for the present work: we provide a systematic construction for generating
such distributions. We will produce a simple scheme for constructing random variables with arbitrarily heavy tails whose minimum has an arbitrarily
light tail.

To formulate our results, we need to introduce some notation and notions. 
Throughout the paper, we use the same notation $F$ for a probability distribution on the real line and for its
distribution function $F(x)$, and we denote by
\[
    \ov{F}(x) = 1 - F(x)
\]
its right tail.
The following definitions are standard and can be found, for example,
in \cite{EKM,FKZ}.

\begin{defn}
    A distribution $F$ (or a random variable $\xi$ following distribution $F$) has
    a \emph{heavy tail} if, for all $\lambda > 0$,
    \[
        \int_{-\infty}^{\infty} e^{\lambda x} F(dx)
        \equiv \mathbb{E} e^{\lambda \xi} = \infty, 
    \]
    and  a \emph{light tail} if 
    $\E e^{\lambda \xi} <\infty $, for some $\lambda >0$.
    
    A distribution $F$ is \emph{long-tailed} if it has a right-unbounded support  
    (i.e. $F(x)<1$ for all $x$) and
    \[
    \lim_{x\to\infty} \frac{\ov{F}(x+1)}{\ov{F}(x)} = 1.
    \]
\end{defn}
It is known (see e.g. \cite{FKZ}) that \emph{any long-tailed distribution is heavy-tailed}.

Here is our first result.

\begin{thm}
\label{theorem:min_of_two}
    Let $F$ be any distribution with right-unbounded support,
    and let $g$ be any non-negative increasing function satisfying
    $\lim_{x \to \infty} g(x) = \infty$.

    Then there exist independent random variables $\xi_1$ and $\xi_2$ such that
    the distribution of $\min(\xi_1, \xi_2)$ is $F$, and
    \[
        \mathbb{E} g(\xi_1) = \mathbb{E} g(\xi_2) = \infty.
    \]
\end{thm}

As a particular corollary, we show the universality of the fact that \emph{any light-tailed random variable may be represented as the minimum of two independent heavy-tailed random variables}.

\begin{cor}\label{Cor1}
Given any light-tailed distribution $F$, one can choose two heavy-tailed distribution $H_1$ and $H_2$ such that, for independent random variables $\xi_1$ with  distribution $H_1$ and $\xi_2$ with distribution $H_2$, the random variable $\eta = \min (\xi_1,\xi_2)$ has distribution $F$.
\end{cor}

Clearly, distributions $H_1$ and $H_2$ in the corollary  cannot be long-tailed (otherwise, the minimum would be also long-tailed and, therefore, heavy-tailed).

Consider the following practical scenario, which illustrates the idea behind the construction of
the two distributions in the proof of this theorem.
Imagine two developer teams, each already overwhelmed with tasks.
A crucial client reports a critical bug that must be fixed without delay.
Initially, the bug is assigned to both teams, but because each is already overloaded,
it lands at the bottom of their respective queues.
Recognizing the urgency, management instructs one team to drop everything and address
the bug immediately, granting them a fixed amount of time.
If that team fails to resolve the issue within the allotted time, the bug is reassigned to
the other team, this time with a more generous time limit, while the first team returns to
its previous work.
If both teams fail in their first attempts, the cycle continues with successively longer time limits
until the bug is finally fixed.

This setup reflects the core idea examined in this paper.
While the expected time for either team to solve the bug might be very large,
suggesting a heavy-tailed distribution, the overall expected resolution time for the bug
remains small, indicating a light-tailed distribution.

We extend this idea to the minimum of several random variables.
To do so, we need a way to compare distributions.
One common way to compare distributions is
via the classical stochastic dominance, see e.g. \cite{MS}:

\begin{defn}
    We say that a distribution $F$ \emph{stochastically dominates} a distribution $G$ if
    \[
        F(x) \leqslant G(x)
        \quad\text{for all } x \in \mathbb{R}.
    \]
\end{defn}

This allows us to establish a more general result concerning the minimum of
many random variables.

\begin{thm}
\label{theorem:min_of_many}
    Let $F$ be any distribution with right-unbounded support,
    and let $g$ be a positive, strictly increasing function satisfying
    $\lim_{x \to \infty} g(x) = \infty$.

    Then, for any $n \in \mathbb{N}$ and any $1 < k \leqslant n$,
    there exist independent
    random variables $\xi_1, \ldots, \xi_n$ such that the following properties hold:

    \begin{enumerate}
        \item For any set of indices 
            \[
                1 \leqslant i_1 < i_2 < \dots < i_k \leqslant n,
            \]
            the distribution $F$ stochastically dominates the distribution of
            $\min(\xi_{i_1}, \dots, \xi_{i_k})$.
        \item For any set of indices 
            \[
                1 \leqslant i_1 < i_2 < \dots < i_{k-1} \leqslant n,
            \]
            we have
            \[
                \E g\Bigl(\min(\xi_{i_1}, \dots, \xi_{i_{k-1}})\Bigr) = \infty.
            \]
    \end{enumerate}
\end{thm}

\section{Proof of Theorem \ref{theorem:min_of_two} and of Corollary \ref{Cor1}, and examples}

\begin{rem}
Without loss of generality, we may reduce our consideration to non-negative random variables. Indeed, for any function $g$ satisfying the conditions of the theorem, $\E g(\xi)$ is finite if and only if $\E g(\max (\xi,0))$ is finite, and, clearly, $\min (\xi_1,\xi_2)^+= \min (\xi_1^+,\xi_2^+)$.
\end{rem}

    We need the following notion.
    The \emph{cumulative hazard function} (for brevity, \emph{cumulative hazard}) of a distribution
    $F$ is defined as
    \[
        R_F(x) = -\log \ov{F}(x).
    \]
The cumulative hazard offers a useful way to quantify how rapidly the tail
$\ov{F}(x)$ decays as $x$ increases.
The remarks below explain the relationship between the cumulative hazard,
a general non-decreasing function, and the minimum of independent random variables.

\begin{rem}
    Any non-decreasing and right-continuous function $R$ that satisfies
    \[
        \lim_{x \rightarrow -\infty} R(x) = 0
        \quad\text{and}\quad
        \lim_{x \rightarrow +\infty} R(x) = \infty
    \]
    is the cumulative hazard function (cumulative hazard) of some distribution $F$, where
    \[
        \ov{F}(x) = e^{-R(x)}.
    \]
\end{rem}

\begin{rem}
    The cumulative hazard of the minimum of two independent random variables $\xi_1$ and $\xi_2$
    equals the sum of their individual cumulative hazards.
    That is, if $\eta = \min(\xi_1, \xi_2)$, then
    \[
        R_{\eta}(x) = R_{\xi_1}(x) + R_{\xi_2}(x).
    \]
\end{rem}

\begin{proof}[Proof of 
Theorem \ref{theorem:min_of_two}]
    Let $H_1$ and $H_2$ be the distribution functions of $\xi_1$ and $\xi_2$.
    The idea of the proof is to let them essentially ``take turns'' on successive
    intervals.
    We divide the positive real axis into consecutive intervals $(a_i, a_{i+1}]$.
    On even-indexed intervals, we ``freeze'' the cumulative hazard $R_1$
    while letting $R_2$ increase together with $R_F$.
    On odd-indexed intervals, we do the opposite: keep $R_2$ fixed and
    let $R_1$ follow $R_F$.
    By carefully stitching these pieces together, we ensure that the sum $R_1 + R_2$
    (i.e., the cumulative hazard of the minimum) matches $R_F$ at all points.
    Hence, $\min(\xi_1, \xi_2)$ has exactly the distribution $F$.
    At each step, we also choose the lengths of these intervals in such a way that each
    distribution's tail is made ``heavy enough.''
    This ``alternating'' approach guarantees that while each random variable is heavy-tailed
    (in the sense prescribed by $g$), their minimum still reproduces the desired lighter-tailed
    distribution $F$.

    We construct the sequence $\{a_i\}_{i=0}^\infty$ and the cumulative hazards
    $R_1, R_2$ inductively. 
    First, set $a_0 = 0$ and define $R_1(x) = R_2(x) = 0$ for all $x \leqslant 0$. 
    Assume $a_l$ is defined.
    We then choose $a_{l+1}$ and specify $R_1$ and $R_2$ on $(a_l, a_{l+1}]$ as follows:
    \begin{itemize}
        \item For even values of $l$, we let
        \[
            R_1(x) = R_1(a_l), \; R_2(x) = R_2(a_l) + \bigl(R_F(x) - R_F(a_l)\bigr)
        \]
        for $x \in (a_l, a_{l+1}]$.
        
        \item If $l$ is odd, we let 
        \[
            R_1(x) = R_1(a_l) + \bigl(R_F(x) - R_F(a_l)\bigr), \; R_2(x) = R_2(a_l)
        \]
        for $x \in (a_l, a_{l+1}]$.
    \end{itemize}
    This ensures $R_1$ and $R_2$ remain right-continuous at each $a_l$.
    We can choose $a_{l+1}$ arbitrarily provided:
    \begin{itemize}
        \item $a_{l+1} - a_l \geqslant 1$, so the intervals $(a_i, a_{i+1}]$ cover all of
        $\mathbb{R}^+$, and, moreover,  
        \item The following inequalities hold:
        \[
            (a_{l+1} - a_l)\exp\{-R_1(g^{-1}(a_l))\} \geqslant 1
                \quad (\text{if $l$ is even}),
        \]
        \[
            (a_{l+1} - a_l)\exp\{-R_2(g^{-1}(a_l))\} \geqslant 1
                \quad (\text{if $l$ is odd}),
        \]
        ensuring the expectation of each constructed random variable is infinite.
    \end{itemize}
    
    Let $H_1$ and $H_2$ be the distributions corresponding to $R_1$ and $R_2$. We verify that these distributions satisfy the stated properties.
    
    First, we show that \emph{the minimum of the random variables has the correct distribution.}
    Let $R_{\min}$ be the cumulative hazard of $\min(\xi_1, \xi_2)$.
    Clearly, $R_{\min}(a_0) = 0 = R_F(a_0)$. 
    
    Assume there is some $l$ such that $R_{\min}(a_l) = R_F(a_l)$.    
    If $l$ is odd, then, for any $x \in (a_l, a_{l+1}]$, we have
    \[
        R_{\min}(x) - R_{\min}(a_l)
        = \bigl(R_1(x) - R_1(a_l)\bigr) + \bigl(R_2(x) - R_2(a_l)\bigr)
    \]
    \[
        = R_F(x) - R_F(a_l),
    \]
    so $R_{\min}(x) = R_F(x)$ in the interval.

    A similar argument shows that $R_{\min}(x) = R_F(x)$ for $x \in (a_l, a_{l+1}]$.
    Hence, $R_{\min} \equiv R_F$ and $\min(\xi_1, \xi_2)$ has distribution $F$.
    
    Secondly, we show that \emph{the distributions of $\xi_1$ and $\xi_2$ are sufficiently heavy.}
    We show this for $\xi_1$. An identical argument applies to $\xi_2$. Observe that
    \[
        \mathbb{E}g(\xi_1)
        = \int_0^\infty \Prob\bigl(g(\xi_1) > t\bigr)\,dt
        = \sum_{j=0}^\infty \int_{a_j}^{a_{j+1}} \Prob\bigl(g(\xi_1) > t\bigr)\,dt
    \]
    \[
        \geqslant \sum_{j=0}^\infty \int_{a_{2j}}^{a_{2j+1}} \Prob\bigl(g(\xi_1) > t\bigr)\,dt.
    \]
    Since $R_1$ is constant on each interval $(a_{2j}, a_{2j+1}]$, we have
    \[
        \int_{a_{2j}}^{a_{2j+1}} \Prob\bigl(g(\xi_1) > t\bigr)\,dt
        = (a_{2j+1} - a_{2j})
        \exp\{-R_1\bigl(g^{-1}(a_{2j})\bigr)\} 
        \geqslant 1.
    \]
    Thus,
    \[
        \mathbb{E}g(\xi_1)
        \geqslant \sum_{j=0}^\infty 1 
        = \infty.
    \]
\end{proof}

\begin{proof}[Proof of Corollary \ref{Cor1}]
One can apply the theorem by taking a light-tailed distribution $F$ and any function $g$ such that $e^{-cx} = o(g(x))$ as $x\to\infty$, for any $c>0$. For example we may take $g(x) =x^+$. Then both $\xi_1$ and $\xi_2$ have infinite first moments and, therefore, all their exponential moments are infinite, too.

For Corollary \ref{Cor1}, there is an alternative argument.
Theorem 2.6 from \cite{FKZ} says that a distribution $F$ is heavy-tailed if and only if its
cumulative hazard satisfies $\liminf_{x\to\infty} R(x)/x= 0$.
Then we can simplify the above proof of the corollary, by defining a sequence $\{a_l\}$ in such
a way that $R_1(a_{2i})/a_{2i+1} \le 1/i$ and $R_2(a_{2i+1})/a_{2i+2}\le 1/i$, for all $i$.

\end{proof}

Let us provide a few examples illustrating how quickly the sequence
$\{a_i\}_{i=0}^\infty$ constructed in the theorem can grow.

\begin{figure}[t]
  \centering
  \subfloat[Example~\ref{example:exponential}. Log-scale plot of the ``minimal'' recursive endpoints $\{a_k^*\}$ (dashed lines) versus our explicit $\{a_k\}$ from the example (solid curve) for the exponential distribution $\overline{F}(x)=e^{-\alpha x}$ with $g(x)=e^{\beta x}$ ($\beta<\alpha$).]{
    \includegraphics[width=0.45\textwidth]{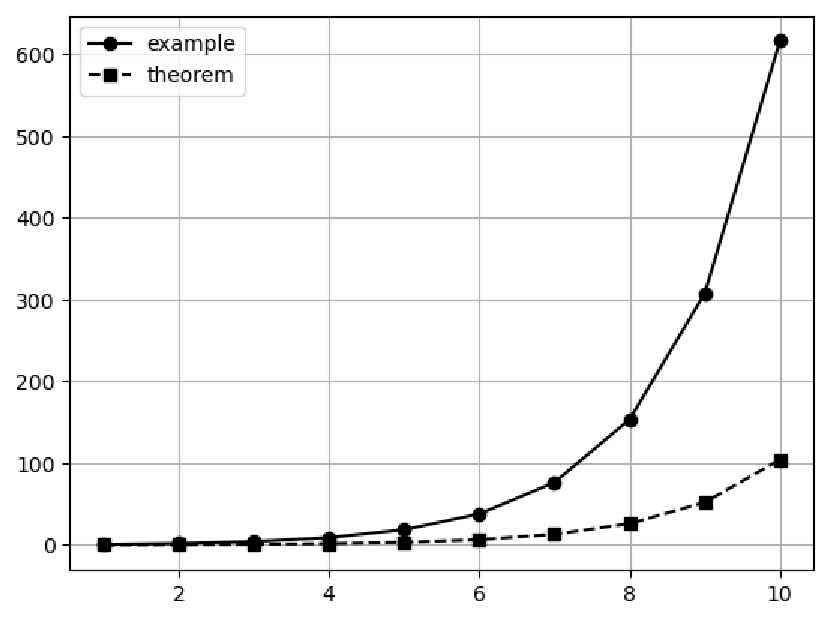}
    \label{ex21}
  }
  \hfill
  \subfloat[Example~\ref{example:power}. Log-scale plot of the ``minimal'' recursive endpoints $\{a_k^*\}$ (dashed lines) versus our explicit $\{a_k\}$ from the example (solid curve) for the power-tailed distribution $\overline{F}(x)=(1+x)^{-\alpha}$ with $g(x)=x^\beta$ ($\beta+1<\alpha$).]{
    \includegraphics[width=0.45\textwidth]{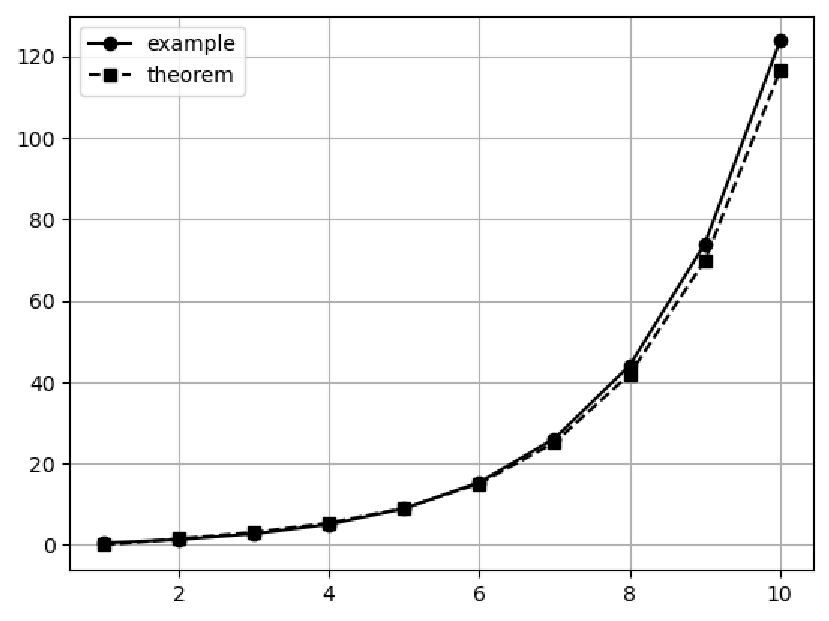}
    \label{ex22}
  }
  \caption{Comparison of explicit $\{a_k\}$ and recursive $\{a_k^*\}$ for two cases: (a) exponential distribution and (b) power-tailed distribution.}
  \label{fig:ex21-ex22}
\end{figure}

To make precise how close our explicit choices of the segment endpoints $a_k$ are to the smallest values allowed by the recursive inequalities in the proof of Theorem \ref{theorem:min_of_two}, we computed a ``minimal'' sequence $\{a_k^*\}$ by enforcing
\[
  a_{k+1}^* - a_k^* 
  = \max\bigl\{1,\exp\!\bigl(R_F\bigl(g^{-1}(a_k^*)\bigr)\bigr)\bigr\},
  \quad a_0^* = 0.
\]
In Figure \ref{fig:ex21-ex22} we then plot both
\[
  \log a_k^*\quad\text{(dashed lines)}
  \quad\text{and}\quad
  \log a_k \quad\text{(solid curves)},
\]
on the same logarithmic (or, where noted, double‐logarithmic) scale. This highlights that, even when our explicit $a_k$ are noticeably larger in some cases, they nevertheless lie on the same asymptotic scale as the theoretical minima.

\begin{ex}
\label{example:exponential}
    Let $F$ be an exponential distribution with $\overline{F}(x) = e^{-\alpha x}$, $x\ge 0$, 
    and let $g(x) = e^{\beta x}$ where $\beta < \alpha$.
    To construct $\xi_1$ and $\xi_2$ according to the theorem, we may take
    \[
        a_k = 2^{\sum_{j=1}^k \left(\frac{\alpha}{\beta}\right)^j}.
    \]
\end{ex}

\begin{ex}
\label{example:power}
    Let $F$ be a distribution with tail $\overline{F}(x) = \frac{1}{(1+x)^\alpha}$, for $x\ge 0$, 
    and let $g(x) = x^\beta$ where $\beta + 1 < \alpha$.
    Then we may also take 
    \[
        a_k = 2^{\sum_{j=1}^k \left(\frac{\alpha}{\beta}\right)^j}.
    \]
\end{ex}

Note, however, that this construction is not always necessary.
Let $F$ be any distribution and define a new distribution $F^\uparrow$ by
\[
    \overline{F}^\uparrow(x) = \sqrt{\overline{F}(x)} \;=\; e^{-R_F(x)/2}.
\]
If $\xi_1$ and $\xi_2$ are independent random variables with distribution $F^\uparrow$,
then $\min(\xi_1,\xi_2)$ will have distribution $F$.
Thus, for the construction in Theorem~\ref{theorem:min_of_two} to be needed,
$\xi_1$ and $\xi_2$ must be heavier-tailed than $F^\uparrow$.

In particular, in Example~\ref{example:exponential}, if $\beta > \alpha/2$,
one can simply take $\xi_1$ and $\xi_2$ to be exponentially distributed with
$\overline{F}(x) = e^{-\alpha x/2}$,
and the minimum of those variables will still have the desired distribution $F$.

\section{Archimedean survival copulas and dependent minima}

In the preceding sections we focused on independent random variables. 
However, this assumption may be considered restrictive in applications.
For example, in actuarial models, risks or lifetimes are often dependent, and copulas provide the standard way to model such dependence.
A particularly convenient class is given by Archimedean survival copulas,
which offer tractable families of dependence models widely used not only in
survival analysis, but also in finance and in many other areas of applied probability
(see, e.g., \cite{Nelsen} for a systematic introduction,
including Archimedean survival copulas in particular).

These copulas also fit naturally with our cumulative-hazard-based construction.
If $(X_1,\ldots,X_n)$ has an Archimedean survival copula with generator $\phi$
(non-increasing, $\phi(1)=0$), then for all $x$,
\[
\phi\bigl(\Prob(\min_{i=1,\ldots,n} X_i > x)\bigr)
= \sum_{i=1}^n \phi\bigl(\Prob(X_i > x)\bigr).
\]
Thus the same additive decomposition that underpins the independent case 
appears here as well, only expressed through the transformed function $\phi$. 
This suggests that one can mimic our earlier scheme by requiring heavy-tailed marginals, 
yet still obtain a prescribed light-tailed distribution for the minimum under dependence. 
A full exploration of these possibilities lies beyond the scope of this paper, 
but the structural similarity indicates that our approach is not limited to independence and may extend naturally to certain dependent settings.
In particular, the following result holds.

\begin{thm}
\label{theorem:min_of_two_copulas}
Let $F$ be any distribution with right-unbounded support,
and let $g$ be any non-negative increasing function satisfying
$\lim_{x \to \infty} g(x) = \infty$.

Suppose $h$ is a strictly increasing function that is continuous in a neighborhood of $0$,
with $h(0)=0$ and $h(1)=1$. Then there exist random variables $\xi_1$ and $\xi_2$ such that
\[
    h\bigl(\ov{F}(x)\bigr)
    = h\bigl(\Prob(\xi_1 > x)\bigr)\, h\bigl(\Prob(\xi_2 > x)\bigr),
\]
and
\[
    \E g(\xi_1) = \E g(\xi_2) = \infty.
\]
\end{thm}
\begin{rem}
We do not impose the usual generator conditions that ensure that
$h^{-1}(h(u)h(v))$ is a genuine survival copula.
If one also wants existence of copula-based joint law, additional conditions on $h$
(e.g. convexity of $-\log h$) are required and lie outside our scope.
\end{rem}

The proof can be obtained by adapting the argument of Theorem~\ref{theorem:min_of_two}.
But instead of decomposing cumulative hazards, we decompose survival functions so that
the copula identity is preserved.
More precisely, on each interval $(a_l,a_{l+1}]$ one proceeds as follows:
\begin{itemize}
    \item For even $l$, hold $\ov{F}_{\xi_1}$ constant at the value $\ov{F}_{\xi_1}(a_l)$,
        and let $\ov{F}_{\xi_2}(x)$ vary so that
        \[
            h\bigl(\ov{F}(x)\bigr)
            = h\bigl(\ov{F}_{\xi_1}(a_l)\bigr)\, h\bigl(\ov{F}_{\xi_2}(x)\bigr),
            \qquad x \in (a_l,a_{l+1}].
        \]
    \item For odd $l$, hold $\ov{F}_{\xi_2}$ constant at the value $\ov{F}_{\xi_2}(a_l)$,
        and let $\ov{F}_{\xi_1}(x)$ vary so that
        \[
            h\bigl(\ov{F}(x)\bigr)
            = h\bigl(\ov{F}_{\xi_1}(x)\bigr)\, h\bigl(\ov{F}_{\xi_2}(a_l)\bigr),
            \qquad x \in (a_l,a_{l+1}].
        \]
\end{itemize}

Since $h$ is continuous near $0$, the procedure can always be carried out initially by selecting $a_1$ sufficiently large and putting $\ov{F}_{\xi_1}(x) = 1$ and $\ov{F}_{\xi_2}(x) = \ov{F}(x)$
for $x \in [0, a_1)$.
The subsequent points $a_2,a_3,\ldots$ of the partition can then be chosen large enough to
ensure that both $\xi_1$ and $\xi_2$ are heavy-tailed in the sense required by $g$.

\section{Proof of Theorem \ref{theorem:min_of_many}}

We again restrict our consideration to non-negative random variables.

\begin{proof}
Let $H_1, H_2, \ldots, H_n$ be the distribution functions of
$\xi_1, \ldots, \xi_n$. We present them via their cumulative hazards,
building on the cumulative hazard $R_F$ of $F$.
We partition $\mathbb{R}^+$ into intervals $\{(a_i, a_{i+1}]\}_{i=0}^\infty$ and,
on each interval, we keep the cumulative hazards corresponding to a selected set of
$k-1$ indices constant.
This ensures that the minimum of the associated random variables exhibits a heavy tail.
By cyclically considering all possible subsets of $k-1$ indices, we guarantee that
every such subset yields a heavy-tailed minimum.
Simultaneously, on every interval, each of the remaining $n - (k - 1)$ cumulative hazards is
allowed to grow at the same rate as $R_F$, ensuring that in every interval there is
at least one cumulative hazard that grows quickly enough.

\emph{Construction of intervals and distributions.}
Let $I_0, \ldots, I_{M-1}$ be the collection of all possible subsets of $k-1$ indices from
$\{1, \ldots, n\}$, where $M = \binom{n}{k-1}$.
Extend these subsets into an infinite sequence $\{I_j\}_{j=0}^\infty$ by repeating them
cyclically; that is, identify $I_i$ with $I_j$ whenever $i \equiv j\ (\mathrm{mod}\, M)$.

Next, we inductively define a sequence $\{a_i\}_{i=0}^\infty$ and cumulative hazards
$R_1, \ldots, R_n$.
Set $a_0 = 0$ and define $R_i(x) = 0$ for all $i$ and for all $x \leqslant 0$.
Assuming the functions are defined up to $a_l$,
choose $a_{l+1}$ and define the functions on the interval $(a_l, a_{l+1}]$ as follows:
\begin{itemize}
    \item For any $i \in I_l$, set $R_i(x) = R_i(a_l)$ for all $x \in (a_l, a_{l+1}]$.
    \item For any $i \notin I_l$, define
    \[
        R_i(x) = R_i(a_l) + \Bigl(R_F(x) - R_F(a_l)\Bigr)
    \]
    for all $x \in (a_l, a_{l+1}]$.
\end{itemize}
This construction ensures that the functions $R_1, \ldots, R_n$ remain right-continuous at
$a_l$.
The point $a_{l+1}$ may be chosen arbitrarily provided that:
\begin{itemize}
    \item $a_{l+1} - a_l \geqslant 1$, so that the intervals $(a_i, a_{i+1}]$ are guaranteed to
        cover $\mathbb{R}^+$.
    \item The following inequality holds
    \[
        (a_{l+1} - a_l)\, \exp\Bigl\{-\sum_{i \in I_l} R_i(g^{-1}(a_l))\Bigr\} \geqslant 1,
    \]
    which will guarantee that the expectation of the minimum for any set of $k-1$ indices is
    infinite.
\end{itemize}
Let
$H_1,\ldots,H_n$ be the distributions corresponding to
the cumulative hazards $R_1,\allowbreak\dots,\allowbreak R_n$.
We now verify that these distributions satisfy the stated properties.

We show now that \emph{the minimum of $k$ random variables is sufficiently light.}
Consider any subset of indices $1 \leqslant i_1 < \dots < i_k \leqslant n$ 
with corresponding independent random variables $\xi_{i_1}, \ldots, \xi_{i_k}$.
Let $R_{\min}$ denote the cumulative hazard of their minimum,
$\min(\xi_{i_1}, \dots, \xi_{i_k})$.
By the construction of $\{I_j\}_{j=0}^\infty$, on any interval $(a_l, a_{l+1}]$
there exists at least one index $j_l \in \{i_1, \dots, i_k\}$ with property $j_l \notin I_l$.
Consequently, for any $x \in (a_l, a_{l+1}]$ we have
\[
    R_{\min}(x) - R_{\min}(a_l) =
    \sum_{j\in \{i_1,\ldots,i_k\}} 
    (R_j(x)-R_j(a_l))
    \geqslant R_{j_l}(x) - R_{j_l}(a_l) = R_F(x) - R_F(a_l).
\]
It follows that $R_{\min}(x) \geqslant R_F(x)$ for all $x \geqslant 0$,
so that $F$ stochastically dominates the distribution of $\min(\xi_{i_1}, \dots, \xi_{i_k})$.

It is left to show that \emph{the minimum of $k-1$ random variables is sufficiently heavy.}
Consider any subset of indices $1 \leqslant i_1 < \dots < i_{k-1} \leqslant n$,
and let $\eta = \min(\xi_{i_1}, \dots, \xi_{i_{k-1}})$.
Then,
\[
    \E g(\eta) = \int_0^\infty \Prob\Bigl(g(\eta) > t\Bigr)\,dt
    = \sum_{j=0}^\infty \int_{a_j}^{a_{j+1}} \Prob\Bigl(g(\eta) > t\Bigr)\,dt.
\]
The chosen set of indices corresponds to $I_m$ for some $0 \leqslant m < M$
and we can use the following lower bound
\[
    \E g(\eta) \geqslant \sum_{\substack{j = 0 \\ j \equiv m\,(\mathrm{mod}\, M)}}^\infty
        \int_{a_j}^{a_{j+1}} \Prob\Bigl(g(\eta) > t\Bigr)\,dt.
\]
Moreover, for each such interval $(a_j, a_{j+1}]$, 
\[
    \int_{a_j}^{a_{j+1}} \Prob\Bigl(g(\eta) > t\Bigr)\,dt
    = (a_{j+1} - a_j)\,\exp\Bigl\{-\sum_{i \in I_j} R_i\bigl(g^{-1}(a_j)\bigr)\Bigr\}.
\]
By the choice of $a_{j+1}$, we ensure that
\[
    (a_{j+1} - a_j)\,\exp\Bigl\{-\sum_{i \in I_j} R_i\bigl(g^{-1}(a_j)\bigr)\Bigr\} \geqslant 1.
\]
Therefore,
\[
    \E g(\eta) \geqslant \sum_{\substack{j = 0 \\ j \equiv m\,(\mathrm{mod}\, M)}}^\infty 1
    = \infty.
\]
\end{proof}

\section*{Acknowledgements}
The authors are grateful to the anonymous referees for their constructive comments and suggestions,
which helped to improve the clarity and scope of the paper.
In particular, we thank them for encouraging us to extend the discussion to dependent settings via
Archimedean survival copulas and for their helpful remarks on presentation, references, and notation.

%
%
%


\begin{thebibliography}{5}

\bibitem{LSK1}
{\sc Leipus, R., Siaulys, J. and Konstantinides, D.} (2023). Minimum of heavy-tailed random variables is not heavy tailed. {\em AIMS Mathematics}
{\bf 8 (6)}, 13066--13072.

\bibitem{EKM}
{\sc Embrechts, P., Klüppelberg, C. and Mikosch, T.} (1997). {\em Modelling Extremal Events for Insurance and Finance}. Springer, Berlin.

\bibitem{FKZ}
{\sc Foss, S., Korshunov, D. and Zachary, S.} (2013). {\em An Introduction to Heavy-Tailed and Subexponential Distributions}, 2nd~edn. Springer, New York.

\bibitem{LSK2}
{\sc Leipus, R., Siaulys, J. and Konstantinides, D.} (2023). {\em Closure Properties for Heavy-Tailed and Related Distributions}, Springer.

\bibitem{MS}
{\sc Stoyan, D.} (1983). {\em Comparison Methods for Queues and Other Stochastic Models}. Wiley, Chichester.

\bibitem{MW}
{\sc Matsui, M. and Watanabe, T.} (2025). Convolution closure properties of subexponential densities. {\em Journal of Mathematical Analysis and Applications}, {\bf 545 (1)}, 129--158.

\bibitem{XFW}
{\sc Xu, H., Foss, S. and Wang, Y.} (2015) Convolution and convolution-root properties of long-tailed distributions. {\em Extremes}, {\bf 18(4)}, 605--628.

\bibitem{Nelsen}
{\sc Nelsen, R. B.} (2006). {\em An Introduction to Copulas}, 2nd~edn. Springer, New York.


\end{thebibliography}
\end{document}